\documentclass[a4paper,12pt]{article}

\usepackage{color}
\advance \textwidth -60truemm\relax

\advance\oddsidemargin -1truein\relax

\evensidemargin=\oddsidemargin

\advance\textheight -60truemm\relax
\advance\topmargin -1truein\relax
\unitlength\textwidth
\divide\unitlength by 150\relax

\usepackage{amsmath,amssymb}
\usepackage{bm}
\bmdefine{\aaa}{a}
\bmdefine{\bbb}{b}
\bmdefine{\mmm}{m}
\bmdefine{\ppp}{p}
\bmdefine{\qqq}{q}
\bmdefine{\uuu}{u}
\bmdefine{\vvv}{v}
\bmdefine{\www}{w}
\bmdefine{\eee}{e}
\bmdefine{\xxx}{x}
\bmdefine{\zerovec}{0}

\newcommand{\CCC}{\mathbb{C}}
\newcommand{\RRR}{\mathbb{R}}

\numberwithin{equation}{section}
\newtheorem{thm}[equation]{Theorem}

\newtheorem{lemma}[equation]{Lemma}

\newtheorem{prop}[equation]{Proposition}
\newtheorem{remark}[equation]{Remark}

\newcommand{\bigzerou}{\smash{\lower1.7ex\hbox{\bg 0}}}

\newcommand{\bigastu}{\smash{\lower1.7ex\hbox{\bg *}}}

\begin{document}

\title { A simple estimation of the maximal rank of tensors with two slices by row and column operations, symmetrization and induction}
\author{Toshio Sakata, Toshio Sumi and Mitsuhiro Miyazaki%
\thanks{Kyushu University, Kyoto University of Education and Kyushu University}
}
\maketitle  
\section{Introduction}

The determination of the maximal ranks of a set of a given type of tensors is a basic problem both in theory and application. 
In statistical applications, the maximal rank is related to the number of necessary parameters to be built in a tensor model.  
JaJa \cite{JA'JA'} and Sumi et.~al \cite{Sumi-Miyazaki-Sakata-1} developed an optimal bound theory based on Kronecker canonical form of the pencil of two matrices.  Theory of matrix pencil is explained in several text book, for example, of Gantmacher~\cite{Gantmacher}.  Atkinson and Lloyd\cite{Atkinson-Lloyd},  
Atkinson and Stephens\cite{Atkinson-Stephens} and Sumi et.~al \cite{Sumi-Miyazaki-Sakata-2} treated the maximal rank of tensors with 3 slices of matrices.  In contrast we use an old theorem, which states that any real matrix can be expressed as a product of two real 
symmetric matrices.  Based on this classical theorem (Bosch~\cite{Bosch})  we will show the tight bound by simple row and column operations and symmetrization and mathematical induction. As far as the authors know, the inductive proof of the tight bound $[3n/2]$ for $2 \times n \times n$ tensors, which has been given by several authors based on eigenvalue theories, is the first result in this filed.  It should be note that the inductive proof is shown to have a great difficulty for odd $n$.  We overcame this in this paper.  In Section 2 we list up several proofs for some particular cases, which are very interesting in themselves and became stepstones of our general proof. 
In Section 3 we will give a proof by symmetrization and an inductive proof for the maximal rank of $2 \times n\times n$. 
Finally, in Section 4, we will generalize the proof for the case of  $2 \times m \times n$ tensors.

\section{Estimation by using row and column operation}
 In this section we list up the bounds, which can be obtained simply by appropriate row and column operations, different for each particular cases. These standalone results became our motivation for more simpler proof than one based on eigenvalues.  Here we denote the set of all $2 \times m\times n$ tensors by $T(2,m,n)$ and the maximal rank of tensors in $T(2,m,n)$ is denoted by shortly $r(2,m,n)$. Also we use the notation $r(T)$ for the rank of a particular tensor $T$.  It should be noted that in this section for almost all cases we consider a $ 2 \times m \times n$ tensor as an object with a slice of $m \times  n$ matrices and therefore  all symbol $\bm{a}, \bm{b}$ and $*$  denote a $2$-dimensional vector and $\bm{0}$ denotes the $2$-dimensional zero vector. Exceptional case is Proposition 2.4, where the symbols denote $3$-dimensional vectors.

\subsection{$2 \times 2 \times n $}

\begin{prop} It holds that 
$r(2,2,2)=3$.  
\end{prop}

\begin{proof}
 $T$ is expressed as  
\[T=\left( 
\begin{array}{ccc}
  \bm{a}  &   \bm{b}           \\ 
  \bm{c}  &   \bm{d}                  \\ 
\end{array}
\right)_{.}
 \]
 Clearly it suffices to prove the proposition   $\bm{a}$ and $\bm{b}$ are independent and $\bm{c}$ is not a constant multiple of $\bm{a}$. Then we can express $T$ as 
 \[T=\left( 
\begin{array}{ccc} 
  \bm{a}  &   \bm{b}           \\ 
  \alpha \bm{a}+\beta \bm{b}  &   \gamma \bm{a} +\delta \bm{b}         
\end{array}
\right)_{.}
 \]
By a row operation and constant multiplication to the 2nd row we have 
\[T=\left( 
\begin{array}{ccc} 
  \bm{a}  &   \bm{b}           \\ 
  \bm{b}  &   \gamma \bm{a} +\delta \bm{b}                  \\ 
\end{array}
\right)_{.}
 \]
If $\delta=0$ we have 
\[T=\left( 
\begin{array}{ccc} 
  \bm{a}  &   \bm{b}           \\ 
  \bm{b}  &   \gamma \bm{a}     \\ 
\end{array}
\right)_{.}
 \]
 and this is decomposed as 
 \[T=\left( 
\begin{array}{ccc} 
  \bm{a}-\bm{b}  &   \bm{b}           \\ 
  \bm{b}  &   \gamma \bm{a}-\bm{b}     \\ 
\end{array}
\right) 
+
\left( 
\begin{array}{ccc} 
  \bm{b}  &   \bm{b}           \\ 
  \bm{b}  &   \bm{b}     \\ 
\end{array}
\right) 
 \]
 and $r(T) \leq 3$
If $\delta \neq 0$, by constant multiplications, we have 
\[T=\left( 
\begin{array}{ccc} 
 \delta \bm{a}  &   \bm{b}           \\ 
  \bm{b}  &   \frac{\gamma}{\delta}  \bm{a}+\bm{b}     \\ 
\end{array}
\right)
 \]
and this is decomposed as  
 \[T=\left( 
\begin{array}{ccc} 
  \delta \bm{a}-\bm{b}  &   \bm{0}           \\ 
  \bm{0}  &     \frac{\gamma}{\delta}  \bm{a}    \\ 
\end{array}
\right) 
+
\left( 
\begin{array}{ccc} 
  \bm{b}  &   \bm{b}           \\ 
  \bm{b}  &   \bm{b}     \\ 
\end{array}
\right)_{.}
 \]
Thus $r(T) \leq 3$. These complete the proof. 
\end{proof}
The next result is somewhat surprising,  because the maximal rank of $T(2,2,3)$ is the same with one of $T(2,2,2)$, nevertheless $T(2,2,3)$ is truly larger than
$T(2,2,2)$.  
\begin{prop} 
$r(2,2,3)=3$
\end{prop}

\begin{proof}
If $T$ is 
\[
\left( \begin{array}{ccc} 
  \bm{0}  &        \bm{0}   &   \bm{0}      \\ 
 *   & *           &   *             \\ 
\end{array} \right)_{,}
 \]
clearly $r(T)=3$. 
So, we assume that T is 
\[
\left( \begin{array}{ccc} 
  \bm{a}  &        \bm{b}   &   \bm{c}      \\ 
 *   & *           &   *             \\ 
\end{array} \right)_{,}
 \]
where $\bm{a} \neq \bm{0}$.  If  both $\bm{b}$ and $\bm{c}$ are multiple of $\bm{a}$, by operation of columns, $T$ becomes
 \[
\left( \begin{array}{ccc} 
  \bm{a}  &      \bm{0}   &   \bm{0}      \\ 
 *   &     \bm{d}           &   \bm{e}             \\ 
\end{array} \right)_{.}
 \]
Then, if  $\bm{d}$ and $\bm{e}$ is independent, by column operation, $T$ becomes 
  \[
\left( \begin{array}{ccc} 
  \bm{a}  &      \bm{0}   &   \bm{0}      \\ 
 *   &     \bm{d}           &   \bm{e}             \\ 
\end{array} \right)
 \]
 and the rank of $T$ is 3.
If $\bm{d}$ and $\bm{e}$ is dependent, by column operation,  T becomes 
  \[
\left( \begin{array}{ccc} 
  \bm{a}  &      \bm{0}   &   \bm{0}      \\ 
 *   &     \bm{0}           &   \bm{e}             \\ 
\end{array} \right)
 \]
or 
  \[
\left( \begin{array}{ccc} 
  \bm{a}  &      \bm{0}   &   \bm{0}      \\ 
 *   &     \bm{d}           &   \bm{0}             \\ 
\end{array} \right)
 \]
 and the rank is 3 in any way.
Next we consider that 
 T is 
\[
\left( \begin{array}{ccc} 
  \bm{a}  &  \bm{b}   &   \bm{0}      \\ 
 *   & *    &  \alpha \bm{a}+\beta \bm{b}             \\ 
\end{array} \right)_{,}
 \]
 where  $\bm{a}$ and $\bm{b}$ are linearly independent. Then, if $\alpha=0$ and $\beta \neq 0$, $T$ becomes by column operations 
 \[
\left( \begin{array}{ccc} 
  \bm{a}  &  \bm{b}   &   \bm{0}      \\ 
 *   & *    &  \bm{b}             \\ 
\end{array} \right)_{.}
 \]
 And further, by column operations, $T$ becomes 
  
\[ \left( \begin{array}{ccc} 
 \bm{a}  &  \bm{b}   &   \bm{0}      \\ 
  \gamma \bm{a}   & \delta \bm{a}    &  \bm{b}   \\     
  \end{array} \right)_{.}
 \]
 If $\gamma=0$, $T$ becomes 
\[ \left( \begin{array}{ccc}  
\bm{a}  &  \bm{b}   &   \bm{0}      \\ 
  \bm{0}   & \delta \bm{a}    &  \bm{b}             \\ 
  \end{array} \right)_{.}
 \]
If $\delta=0$, the rank is 3 and we assume that $\delta \neq 0$. Then, multiplications by constants to the 2nd rows and the 2nd column,  $T$ becomes 
 by column operations, $T$ becomes 
\[  \left( \begin{array}{ccc}  
 \bm{a}  &  \bm{b}   &   \bm{0}      \\ 
   \bm{0}   &  \bm{a}    &  \bm{b}             \\ 
   \end{array} \right)_{.}
 \]
 Adding 1st column and 3rd column to 2nd column, $T$ becomes
 \[  \left( \begin{array}{ccc} 
   \bm{a}  &  \bm{a}+\bm{b}   &   \bm{0}      \\ 
   \bm{0}   &  \bm{a}+\bm{b}    &  \bm{b}             \\ 
   \end{array} \right)
 \]

 and the rank of $T$ is 3. If $\alpha=0$ and $\beta=0$, $T$ becomes $ 2 \times 2 \time 2$ and of rank 3.
For the case of $\gamma \neq 0$ and $\delta=0$, a similar argument proves that the rank of $T$ is 3.
If $\gamma \neq 0 $ and $\delta \neq 0$ in  
 \[ T=
\left( \begin{array}{ccc} 
  \bm{a}  &  b   &   \bm{0}      \\ 
 *   & *    &  \alpha \bm{a}+\beta \bm{b}             \\ 
\end{array} \right)_{,}
 \]
 by column operations, $T$ becomes
  \[ T=
\left( \begin{array}{ccc} 
  \bm{a}  &  \bm{b}   &   \bm{0}      \\ 
  \gamma \bm{a}   & \delta \bm{b}    &  \alpha \bm{a}+\beta \bm{b}             \\ 
\end{array} \right)_{,}
 \]
 which is clearly of rank 3.  These completes the proof of the proposition.
\end{proof}

\begin{prop}
$r(2,2,p)=4 for  q \ge  4$.
\end{prop}
\begin{proof} 
The proof  of their fact is easy and omitted.
\end{proof} 

\subsection{$2 \times 3 \times n$}
First we show $r(2,3,3) \leq 4$. 
\begin{prop}
$r(2, 3, 3) \leq 4$. 
\end{prop}

\begin{proof}  Here we use the symmetrization method. We can assume that $T$ is \[T=
 \left( \begin{array}{ccc} 
 \bm{a}_{1}  & * & *     \\ 
 *      & *             &              *               \\ 
 *     &  *         &   *   \\ 
\end{array} \right)_{,}
\]
where $\bm{a}_{1} \neq \bm{0}$.  If all vectors in the first row are constant multiples of $\bm{a}_{1}$, by column operations, $T$ becomes
\[T=
 \left( \begin{array}{ccc} 
 \bm{a}_{1}  & \bm{0} & \bm{0}     \\ 
 *      & *             &              *               \\ 
 *     &  *         &   *   \\ 
\end{array} \right)
\]
and then $r(T) \leq r(2,2,3)+1=4$.  Hence we can assume that $T$ is   
\[T=
 \left( \begin{array}{ccc} 
 \bm{a}_{1}  & \bm{b}_{1}  & \bm{0}     \\ 
 *      & *             &              *               \\ 
 *     &  *         &   *   \\ 
\end{array} \right)_{,}
\]
where $\bm{a}_{1},\bm{b}_{1}$ are linearly independent, where (1,3) cell becomes $\bm{0}$ by column operations. By the same argument $T$ becomes 
\[T=
 \left( \begin{array}{ccc} 
 \bm{a}_{1}  & \bm{b}_{1}  &  \bm{0}    \\ 
 \bm{b}_{1}      & *             &              *               \\ 
 \bm{0}     &  *         &   \bm{a}_{2}   \\ 
\end{array} \right)_{,}
\]
where $\bm{b}_{1} $ in (2,1) cell  and (1,2) cell can be taken  identical vectors by constant multiplications.  
If $\bm{a}_{2}=\bm{0}$, then $r(T) \leq r(2,2,3)+1=4$, and so we assume that  $\bm{a}_{2} \neq \bm{0}$.
Then by column operation, $T$ becomes 
\[T=
 \left( \begin{array}{ccc} 
 \bm{a}_{1}  & \bm{b}_{1}  & \bm{0}    \\ 
 \bm{b}_{1}      & *             &              \alpha \bm{b}_{2}               \\ 
 \bm{0}     &   \beta \bm{b}_{2}         &   \bm{a}_{2}   \\ 
\end{array} \right)_{,}
\]
where $\bm{b}_{2}$ is perpendicular to $\bm{a}_{2}$.  Since $\alpha\beta =0$ can be excluded,  by multiplying $1/\beta$ to the 3rd row and multiplying  $1/\alpha$ to the 3rd column, $T$ becomes
\[T=
 \left( \begin{array}{ccc} 
 \bm{a}_{1}  & \bm{b}_{1}  & \bm{0}    \\ 
 \bm{b}_{1}      & *             &       \bm{b}_{2}               \\ 
 \bm{0}     &  \bm{b}_{2}         &   \bm{a}'_{2}   \\ 
\end{array} \right)_{,}
\]
 which is symmetric. First diagonalizing the lower matrix by an orthogonal matrix, and after  multiplying $-1$ if necessary, if  adding a vector in  a diagonal cell, 
the lower matrix can be positive diagonal matrix and therefore can be the identity matrix by a diagonal multiplication of a positive diagonal matrix from left and right transformation. For this operations the upper matrix remains symmetric and so by multiplying an orthogonal matrix to the both matrix we have a diagonal matrix simultaneously on the upper and lower matrices.  Therefore the rank is 3, and after  deleting the added diagonal tensor, the rank of tensor is 4.
\end{proof}

\begin{prop}
$r(2,3,4) \leq 5$.
\end{prop}

\begin{proof}
We can start by 
\[T=
 \left( \begin{array}{cccc} 
 \bm{a}  & \bm{b}  &  \bm{0} &  \bm{0}    \\ 
 \bm{b}  & *  &   * & *             \\ 
 \bm{0}     & * & * &   *   \\ 
\end{array} \right)_{,}
\]
where $\bm{a}$ and $\bm{b}$ are independent. Then by row and column operations, $T$ becomes 
\[T=
 \left( \begin{array}{cccc} 
 \bm{a}  & \bm{b}  & \bm{b} & \bm{b}    \\ 
 \bm{b}  & *  & * & *             \\ 
\bm{ b}  & *  & * & *   \\ 
\end{array} \right)
\]
and decompose this as
\[T=
 \left( \begin{array}{cccc} 
 \bm{a-b}  & \bm{0}  & \bm{0} &\bm{ 0}    \\ 
 \bm{0}  & \bm{0}  & \bm{0} & \bm{0}             \\ 
 \bm{0}  & \bm{0}  & \bm{0} & \bm{0}   \\ 
\end{array} \right)
+
\left( \begin{array}{cccc} 
 \bm{b}  & \bm{b}  & \bm{b} & \bm{b}    \\ 
 \bm{b}  & \bm{b}  & \bm{b} & \bm{b}             \\ 
 \bm{b}  & \bm{b}  & \bm{b} & \bm{b}   \\ 
\end{array} \right)
+
\left( \begin{array}{cccc} 
 \bm{0}  & \bm{0}  & \bm{0} & \bm{0}    \\ 
 \bm{0}  & *-\bm{b}  & *-\bm{b} & *-\bm{b}             \\ 
 \bm{0}  & *-\bm{b}  & *-\bm{b} & *-\bm{b}   \\ 
\end{array} \right)
\]
and from this, we have the estimate,  
$$
r(2,3,4) \leq 1+1+r(2,2,3)=2+3=5
$$
\end{proof}

\begin{prop}
$r(2,3,5)  \leq 5$
\end{prop}
\begin{proof}
Here exceptionally we consider the tensor as a object with three slices of $2 \times 5$ matrices. Thus each symbol denotes a $3$-dimensional vector. 

If all the vectors of the first row are dependent, by column operations, 
\[T=
\left( \begin{array}{ccccc} 
          \bm{a}      & \bm{0}  & \bm{0}    &  \bm{0}  &\bm{ 0}        \\ 
          *    &    \bm{b}   & \bm{c}    &   \bm{d}          &  \bm{e}              \\ 
\end{array} \right)_{.}
\]
Then, we have the estimate of $1+r(1,3,5)=1+3=4$.

Next if  the vector space spanned by the vectors in the first row is $2$-dimensional, by column operations, $T$ becomes
\[T=
\left( \begin{array}{ccccc} 
          \bm{a}      & \bm{b}  & \bm{0}    &  \bm{0}  & \bm{0}        \\ 
          *    &   *   & \bm{c}    &  \bm{ d}          &  \bm{e}              \\ 
\end{array} \right)_{.}
\]
If $\dim\langle c,d,e \rangle<3$, the case reduces to the case of $2\times 3\times 4$ and by Proposition 2.3 the maximal rank  is estimated as
5.

If the vector space $\langle \bm{c},\bm{d},\bm{e}\rangle$ is $3$-dimensional, by column operations, $T$ becomes 
 \[T=
\left( \begin{array}{ccccc} 
          \bm{a}      & \bm{b}  & \bm{0}    &  \bm{0}  & \bm{0}        \\ 
          \bm{0}    &    \bm{0}   & \bm{c}    &   \bm{d}    &  \bm{e}              \\ 
\end{array} \right)_{,} 
\]  
and the rank of $T$ is at most $5$.  Finally, the remaining case is one where  both the vector spaces generated vectors in the first row and in the second row are
$3$-dimensional. Then, by column operations,  $T$ becomes 
   \[T=
\left( \begin{array}{ccccc} 
          \bm{a}      & \bm{b}  & \bm{c}    &  \bm{0}  & \bm{0}        \\ 
          \bm{d}    &    \bm{e}   &  \bm{f}   &   \bm{g}    &  \bm{h}              \\ 
\end{array} \right)_{.}
\]   
%%%%%
If $\bm{g}$ and $\bm{h}$ are dependent, by column operations, $T$ becomes 
   \[T=
\left( \begin{array}{ccccc} 
          \bm{a}      & \bm{b}  &\bm{ c}    &  \bm{0}  & \bm{0}        \\ 
          \bm{d}    &    \bm{e}   & \bm{f}  &   \bm{g}    &  \bm{0}              \\ 
\end{array} \right)_{.}
\]   
and $T$ can be viewed as $2 \times 3 \times 4$ and the rank is at most $5$. So we assume $\bm{g}$ and $\bm{h}$ are independent. Since $\bm{a}, \bm{b}$ and $\bm{c}$ 
are assumed independent, by column operations, $T$ becomes 
   \[T=
\left( \begin{array}{ccccc} 
          \bm{a}'      & \bm{b}'  & \bm{c}'    &  \bm{0}  & \bm{0}        \\ 
          \bm{0}    &    \bm{0}   &  \bm{f} ' &   \bm{g}    &  \bm{h}              \\ 
\end{array} \right)_{,}
\]   
where $ \bm{a}', \bm{b}'$ and  $\bm{c}'$ are independent and   $ \bm{f}', \bm{g}$ and  $\bm{h}$ are independent. 
Then there is a vector $\bm{z}$ such that
\[T=
\left( \begin{array}{ccccc} 
          \bm{a}'      & \bm{b}'  & \alpha_{1} \bm{a}' +\beta_{1} \bm{b}' +\gamma_{1} \bm{z}   &  \bm{0}  & \bm{0}        \\ 
          \bm{0}    &    \bm{0}   &   \alpha_{2} \bm{a}' +\beta_{2} \bm{b}' +\gamma_{2} \bm{z}  &   \bm{g}    & \bm{ h}              \\ 
\end{array} \right)
\]   
for suitable $\alpha_1,\alpha_2,\beta_1,\beta_2,\gamma_1,\gamma_2$.
Hence, by column operations, $T$ becomes
\[T=
\left( \begin{array}{ccccc} 
          \bm{a}'      & \bm{b}'  & \gamma_{1} \bm{z}   &  \bm{0}  & \bm{0}        \\ 
          \bm{0}    &    \bm{0}   & \gamma_{2} \bm{z}  &   \bm{g}    &  \bm{h}              \\ 
\end{array} \right)_{.}
\]   
Thus the rank of $T$ is at most $5$. This completes the proof.  
\end{proof}

\subsection{$ 2 \times 4 \times 4$}
\begin{prop} \label{Pro244}
$r(2,4,4) \leq 6$.
\end{prop}

\begin{proof}
We start from 
\[T=
\left( \begin{array}{cccc} 
          \bm{a}      & \bm{b}  & *    &  *         \\ 
          \bm{b}    &    *        &   *          &  *              \\ 
          *    &   *         &   *          &   *      \\ 
          *    &   *         &   *          &   *                 \\ 
\end{array} \right)_{,}
\]
where $\bm{a}$ and $\bm{b}$ are linearly independent, because otherwise the 1st row or the 1st column has the form 
of $(\bm{a},\bm{0},\bm{0},\bm{0})$ and the tensor $T$ can be decomposed as the sum of  a element of $T(2,3,4)$ and a element of $T(1,1,4)$ and $r(T) \leq 5+1=6$.   
In this form, by column operation and row operations, $T$ becomes
\[T=
\left( \begin{array}{cccc} 
          \bm{a}    & \bm{b}  &  \bm{0}   &  \bm{0}         \\ 
          \bm{b}    & *  &  *   &  *                \\ 
          \bm{0}    & *  &  *   &  *      \\ 
          \bm{0}    & *  &  *   &  *                 \\ 
\end{array} \right)
\]
By adding the 2nd row (resp. column) to the 3rd row (resp. column) and the 4th row (resp. column), $T$ becomes 
\[T=
\left( \begin{array}{cccc} 
          \bm{a}    &  \bm{b}  &   \bm{b}   &   \bm{b}         \\ 
          \bm{b}    &  *  &   *   &   *         \\ 
          \bm{b}    &  *  &   *   &   *         \\ 
          \bm{b}    &  *  &   *   &   *         \\ 
\end{array} \right)_{.}
\]
Then we decompose $T$ as
\[
\left( \begin{array}{cccc} 
          \bm{a}-\bm{b}  &   \bm{0}   &  \bm{0}   &  \bm{0}         \\ 
          \bm{0}    &   \bm{0}   &  \bm{0}   &  \bm{0}         \\ 
          \bm{0}    &   \bm{0}   &  \bm{0}   &  \bm{0}         \\ 
          \bm{0}    &   \bm{0}   &  \bm{0}   &  \bm{0}         \\ 
\end{array} \right)
+
\left( \begin{array}{cccc} 
          \bm{b}    &   \bm{b}   &  \bm{b}  &  \bm{b}         \\ 
          \bm{b}    &   \bm{b}   &  \bm{b}  &  \bm{b}         \\ 
          \bm{b}    &   \bm{b}   &  \bm{b}  &  \bm{b}        \\ 
          \bm{b}    &   \bm{b}   &  \bm{b}  &  \bm{b}         \\ 
\end{array} \right)
+
\left( \begin{array}{cccc} 
          \bm{0}   &  \bm{0}    & \bm{0}   &  \bm{0}        \\ 
          \bm{0}   &  *-\bm{b}  & *-\bm{b} &  *-\bm{b}        \\ 
          \bm{0}   &  *-\bm{b}  & *-\bm{b} &  *-\bm{b}     \\ 
          \bm{0}   &  *-\bm{b}  & *-\bm{b} &  *-\bm{b}         \\ 
\end{array} \right)_{.}
\]
From this decomposition we have that $r(T) \leq 1+1+r(2,3,3)=6$. 
\end{proof}

\subsection{$2 \times 5 \times 5$}

\begin{prop}
$r(2,5,5) \leq 7$
\end{prop}
\begin{proof}
Let $T=(A_{1}:A_{2})$. \\
(Case 1.) If $A_{1}$ or $A_{2}$ is non-singular the proof is easy by using symmetrization. For the symmetrization see in the subsection~4.1. \\ 
(Case 2.) If both of $A_{1}$ and $A_{2}$ is singular and $A_{1}$ or $A_{2}$ is of rank less than equal to $3.$ Here we assume that the rank of $A_{2}$ is less than or equal to $3.$ Then by appropriate transformation, $T$ becomes
\[
\left( \begin{array}{ccccc} 
       * & * & * & * & * \\  
       * & * & * & * & * \\
       * & * & * & * & * \\
       * & * & * & * & * \\
       * & * & * & * & * 
       \end{array}
       \right)
       :
       \left(
       \begin{array}{ccccc} 
       * & * & * & 0 & 0 \\  
       * & * & * & 0 & 0 \\
       * & * & * & 0 & 0 \\
       0 & 0 & 0 & 0 & 0 \\
       0 & 0 & 0 & 0 & 0 
       \end{array}
       \right)_{.}
       \]
       Further this is decomposed into 
  \[
\left( \begin{array}{ccccc} 
       * & * & * & 0 & 0 \\  
       * & * & * & 0 & 0 \\
       * & * & * & 0 & 0 \\
       * & * & * & 0 & 0 \\
       * & * & * & 0 & 0 
       \end{array}
       \right)
       :
       \left(
       \begin{array}{ccccc} 
       * & * & * & 0 & 0 \\  
       * & * & * & 0 & 0 \\
       * & * & * & 0 & 0 \\
       0 & 0 & 0 & 0 & 0 \\
       0 & 0 & 0 & 0 & 0
       \end{array}
       \right)
       +
    \left( \begin{array}{ccccc}    
       0 & 0 & 0 & * & * \\  
       0 & 0 & 0 & * & * \\
       0 & 0 & 0 & * & * \\
       0 & 0 & 0 & * & * \\
       0 & 0 & 0 & * & * 
       \end{array}
       \right)
       :
       \left(
       \begin{array}{ccccc} 
       0 & 0 & 0 & 0 & 0 \\  
       0 & 0 & 0 & 0 & 0 \\
       0 & 0 & 0 & 0 & 0 \\
       0 & 0 & 0 & 0 & 0 \\
       0 & 0 & 0 & 0 & 0 
       \end{array}
       \right)_{.}
       \]     
     Hence $r(T) \leq r(2,3,5)+2=5+2=7$. \\
Case(3). Both of $A_{1}$ and $A_{2}$ is of rank $(n-1)$.   \\
We can start from   
\[T=
\left( \begin{array}{ccccc} 
       * & * & * & * & * \\  
       * & * & * & * & * \\
       * & * & * & * & * \\
       * & * & * & * & * \\
       * & * & * & * & x 
       \end{array}
       \right)
       :
       \left(
       \begin{array}{ccccc} 
       1 & 0 & 0 & 0 & 0 \\  
       0 & 1 & 0 & 0 & 0 \\
       0 & 0 & 1 & 0 & 0 \\
       0 & 0 & 0 & 1 & 0 \\
       0 & 0 & 0 & 0 & 0 
       \end{array}
       \right)_{.}
       \]  
       If $x \neq 0$, $T$ is equivalent to   
       \[T=
\left( \begin{array}{ccccc} 
       * & * & * & * & 0 \\  
       * & * & * & * & 0 \\
       * & * & * & * & 0 \\
       * & * & * & * & 0 \\
       0 & 0 & 0 & 0 & 1 
       \end{array}
       \right)
       :
       \left(
       \begin{array}{ccccc} 
       1 & 0 & 0 & 0 & 0 \\  
       0 & 1 & 0 & 0 & 0 \\
       0 & 0 & 1 & 0 & 0 \\
       0 & 0 & 0 & 1 & 0 \\
       0 & 0 & 0 & 0 & 0 
       \end{array}
       \right)_{.}
       \]
From this $r(T) \leq r(2,4,4)+1=6+1=7$. Therefore,we assume $x=0$ and we have 
        \[T=
\left( \begin{array}{ccccc} 
       * & * & * & * & x_{1} \\  
       * & * & * & * & x_{2} \\
       * & * & * & * & x_{3} \\
       * & * & * & * & x_{4} \\
       y_{1} & y_{2} & y_{3} & y_{y} & 0 
       \end{array}
       \right)
       :
       \left(
       \begin{array}{ccccc} 
       1 & 0 & 0 & 0 & 0 \\  
       0 & 1 & 0 & 0 & 0 \\
       0 & 0 & 1 & 0 & 0 \\
       0 & 0 & 0 & 1 & 0 \\
       0 & 0 & 0 & 0 & 0 
       \end{array}
       \right)_{.}
       \] 
       If $(x_{1},x_{2},x_{3},x_{4})=(0,0,0,0)$, $T$ becomes $(2,4,5)$ type and $r(T) \leq r(2,4,4)+1=7$. So, we assume
 that  $(x_{1},x_{2},x_{3},x_{4}) \neq (0,0,0,0).$ Similarly we can assume that $(y_{1},y_{2},y_{3},y_{4}) \neq (0,0,0,0).$   
Then after appropriate transpositions of rows and columns and equivalent transformations and constant multiplications, we have 
\[
\left( \begin{array}{ccccc} 
       * & * & * & 0 & 0 \\  
       * & * & * & 0 & 0 \\
       * & * & * & 0 & 0 \\
       0 & 0 & 0 & 0 & 1 \\
       0 & 0  & 0  & 1 & 0 
       \end{array}
       \right)
       :
       \left(
       \begin{array}{ccccc} 
       * & * & * & * & 0 \\  
       * & * & * & * & 0 \\
       * & * & * & * & 0 \\
       * & * & *& * & 0 \\
       0 & 0 & 0 & 0 & 0 
       \end{array}
       \right)_{.}
       \] 
     Since $A_{1}$ is of rank 4, without loss of generality, we can assume that the 1st and the 2nd column are independent and so 
the 3rd column can be the zero vector by using the 1st and the 2nd columns.  After that,  without loss of generality,  we can assume that the 1st and the 2nd rows are independent and so the 3rd column can be the zero vector, also.  Thus we have 
\[
\left( \begin{array}{ccccc} 
       x_{11} & x_{12} & 0 & 0 & 0 \\  
       x_{21} & x_{22} & 0 & 0 & 0 \\
       0 & 0 & 0 & 0 & 0 \\
       0 & 0 & 0 & 0 & 1 \\
       0 & 0  & 0  & 1 & 0 
       \end{array}
       \right)
       :
       \left(
       \begin{array}{ccccc} 
       * & * & * & * & 0 \\  
       * & * & * & * & 0 \\
       * & * & * & * & 0 \\
       * & * & *& * & 0 \\
       0 & 0 & 0 & 0 & 0 
       \end{array}
       \right)_{,}
       \]  
where  
\[
T=\left(\begin{array}{cc}
       x_{11} & x_{12} \\
       x_{21} & x_{22} \\
\end{array}   
\right)
\]
is non singular. By multiplying the matrix from the left 
\[
\left(\begin{array}{ccc}
       T^{-1} & 0_{22} & 0_{21} \\
       0_{22} & E_{22} & 0_{21} \\
       0_{12} & 0_{12} & 1  
\end{array}   
\right)_{,}
\]    
we reach to the following,
 \[
\left( \begin{array}{ccccc} 
       1 & 0 & 0 & 0 & 0 \\  
       0 & 1 & 0 & 0 & 0 \\
       0 & 0 & 0 & 0 & 0 \\
       0 & 0 & 0 & 0 & 1 \\
       0 & 0  & 0  & 1 & 0 
       \end{array}
       \right)
       :
       \left(
       \begin{array}{ccccc} 
       * & * & * & * & 0 \\  
       * & * & * & * & 0 \\
       * & * & * & * & 0 \\
       * & * & *& * & 0 \\
       0 & 0 & 0 & 0 & 0 
       \end{array}
       \right)_{.}
       \]  
We write this as
   \[
\left( \begin{array}{ccccc} 
       1 & 0 & 0 & 0 & 0 \\  
       0 & 1 & 0 & 0 & 0 \\
       0 & 0 & 0 & 0 & 0 \\
       0 & 0 & 0 & 0 & 1 \\
       0 & 0  & 0  & 1 & 0 
       \end{array}
       \right)
       :
       \left(
       \begin{array}{ccccc} 
       a_{11} & a_{12} & b_{11} & b_{12} & 0 \\  
       a_{21} & a_{22} & b_{21} & b_{22} & 0 \\
       c_{11} & c_{21} & d_{11} & d_{12} & 0 \\
       c_{21} & c_{22} & d_{21}& d_{22} & 0 \\
       0 & 0 & 0 & 0 & 0 
       \end{array}
       \right)_{.}
       \]   
       If $d_{22} \neq 0$, first we decompose as 
   \[\left( \begin{array}{ccccc} 
       0 & 0 & 0 & 0 & 0 \\  
       0 & 0 & 0 & 0 & 0 \\
       0 & 0 & 0 & 0 & 0 \\
       0 & 0 & 0 & 0 & 1 \\
       0 & 0  & 0  & 1 & 0 
       \end{array}
       \right)
       :
       \left(
       \begin{array}{ccccc} 
       0 & 0 & 0 & 0 & 0 \\  
       0 & 0 & 0 & 0 & 0 \\
       0 & 0 & 0 & 0 & 0 \\
       0 & 0 & 0 & 0 & 0 \\
       0 & 0  & 0  & 0 & 0 
       \end{array}
       \right)
       +
       \left( \begin{array}{ccccc} 
       1 & 0 & 0 & 0 & 0 \\  
       0 & 1 & 0 & 0 & 0 \\
       0 & 0 & 0 & 0 & 0 \\
       0 & 0 & 0 & 0 & 0 \\
       0 & 0  & 0  & 0 & 0 
       \end{array}
       \right)
       :
       \left(
       \begin{array}{ccccc} 
       a_{11} & a_{12} & b_{11} & b_{12} & 0 \\  
       a_{21} & a_{22} & b_{21} &_{22} & 0 \\
       c_{11} & c_{21} & d_{11} & d_{12} & 0 \\
       c_{21} & c_{22} & d_{21}& d_{22} & 0 \\
       0 & 0 & 0 & 0 & 0 
       \end{array}
       \right)_{.}
       \]   
Then for the second tensor, by appropriate transformations, we have 
          \[
\left( \begin{array}{ccccc} 
       1 & 0 & 0 & 0 & 0 \\  
       0 & 1 & 0 & 0 & 0 \\
       0 & 0 & 0 & 0 & 0 \\
       0 & 0 & 0 & 0 & 0 \\
       0 & 0  & 0  & 0 & 0 
       \end{array}
       \right)
       :
       \left(
       \begin{array}{ccccc} 
       a_{11} & a_{12} & b_{11} & 0 & 0 \\  
       a_{21} & a_{22} & b_{21} & 0 & 0 \\
       c_{11} & c_{21} & d_{11} & 0 & 0 \\
       0 & 0 &  0& d_{22} & 0 \\
       0 & 0 & 0 & 0 & 0 
       \end{array}
       \right)_{.}
       \]   
      and $r(T) \leq r(2,3,3)+1+2=7.$
Thus we assume that $d_{22}=0$. If $d_{12} \neq 0$, by adding the 3rd row to the 4th row, $d_{22}$ becomes $\neq 0$. Also, if $d_{21} \neq 0$,  adding the 3rd column to the 4th column, $d_{22}$ becomes $\neq 0$. These cases is already excluded, and so we assume  
 that $d_{12}=d_{21}=0$. If $d_{11}\neq 0 $, adding the 3rd column to the 4th column and then adding the 3rd row to the 4th row, we have that $d_{22} \neq 0,$ which is also already excluded. From these argument we can  assume that $d_{11}=d_{12}=d_{21}=d_{22}=0$. So, we have  
\[  \left( \begin{array}{ccccc} 
       1 & 0 & 0 & 0 & 0 \\  
       0 & 1 & 0 & 0 & 0 \\
       0 & 0 & 0 & 0 & 0 \\
       0 & 0 & 0 & 0 & 1 \\
       0 & 0  & 0  & 1 & 0 
       \end{array}
       \right)
       :
       \left(
       \begin{array}{ccccc} 
       a_{11} & a_{12} & b_{11} & b_{12} & 0 \\  
       a_{21} & a_{22} & b_{21} & b_{22} & 0 \\
       c_{11} & c_{21} &  0 & 0 & 0 \\
       c_{21} & c_{22} &  0&  0 & 0 \\
       0 & 0 & 0 & 0& 0 
       \end{array}
       \right)_{.}
       \]  
Since $A_{2}$ is of rank $4$,   
\[B=
\left( \begin{array}{cc}
b_{11} & b_{12} \\
b_{21} & b_{22}
\end{array}
\right)
\ \ and \ \ 
C=\left( \begin{array}{cc}
c_{11} & c_{12} \\
c_{21} & c_{22}
\end{array}
\right)
\]
are both non-singular.  
By multiplying 
 \[
\left( \begin{array}{ccc}
E_{22} & 0 _{22} & 0_{11} \\
B^{-1}A & E_{22} & 0_{21}   \\
0_{12} & 0_{12} & 1
\end{array}
\right)
\]
from the right
we have 
\[  \left( \begin{array}{ccccc} 
       1 & 0 & 0 & 0 & 0 \\  
       0 & 1 & 0 & 0 & 0 \\
       0 & 0 & 0 & 0 & 0 \\
       0 & 0 & 0 & 0 & 1 \\
       * & *  & 0  & 1 & 0 
       \end{array}
       \right)
       :
       \left(
       \begin{array}{ccccc} 
         0 & 0 & b_{11} & b_{12} & 0 \\  
         0 & 0 & b_{21} & b_{22} & 0 \\
       c_{11} & c_{21} &  0 & 0 & 0 \\
       c_{21} & c_{22} &  0&  0 & 0 \\
       0 & 0 & 0 & 0& 0 
       \end{array}
       \right)_{.}
       \]  
 
By multiplying 
 \[
\left( \begin{array}{ccc}
E_{22} & 0 _{22} & 0_{11} \\
0_{21} & B^{-1} & 0_{21}   \\
0_{12} & 0_{12} & 1
\end{array}
\right)
\]
   and 
\[ 
\left( \begin{array}{ccc}
E_{22} & 0 _{22} & 0_{11} \\
0_{21} & C^{-1} & 0_{21}   \\
0_{12} & 0_{12} & 1
\end{array}
\right)
\]
from right and left respectively, we have
\[  \left( \begin{array}{ccccc} 
       1 & 0 & 0 & 0 & 0 \\  
       0 & 1 & 0 & 0 & 0 \\
       0 & 0 & 0 & 0 & c \\
       0 & 0 & 0 & 0 & d \\
       * & *  & *  & * & 0 
       \end{array}
       \right)
       :
       \left(
       \begin{array}{ccccc} 
       0 & 0 & 1 & 0 & 0 \\  
       0 & 0 & 0  & 1 & 0 \\
       1 & 0 &  0 & 0 & 0 \\
       0 & 1 &  0&  0 & 0 \\
       0 & 0 & 0 & 0& 0 
       \end{array}
       \right)_{.}
       \]  
By column changes we have finally, 
\[ T= \left( \begin{array}{ccccc} 
       0 & 0 & 1 & 0 & 0 \\  
       0 & 0 & 0 & 1 & 0 \\
       0 & 0 & 0 & 0 & c \\
       0 & 0 & 0 & 0 & d \\
       * & *  & *  & * & 0 
       \end{array}
       \right)
       :
       \left(
       \begin{array}{ccccc} 
        1 & 0 & 0 & 0 & 0 \\  
        0 & 1 & 0  & 0 & 0 \\
       0 & 0 &  1 & 0 & 0 \\
       0 & 0 &  0&  1 & 0 \\
       0 & 0 & 0 & 0& 0 
       \end{array}
       \right)_{.}
       \]
If $c \neq 0,$ multiplication of a constant to the 5th column, we have
\[ T= \left( \begin{array}{ccccc} 
       0 & 0 & 1 & 0 & 0 \\  
       0 & 0 & 0 & 1 & 0 \\
       0 & 0 & 0 & 0 & 1 \\
       0 & 0 & 0 & 0 & d \\
       * & *  & *  & * & 0 
       \end{array}
       \right)
       :
       \left(
       \begin{array}{ccccc} 
        1 & 0 & 0 & 0 & 0 \\  
        0 & 1 & 0  & 0 & 0 \\
       0 & 0 &  1 & 0 & 0 \\
       0 & 0 &  0&  1 & 0 \\
       0 & 0 & 0 & 0& 0 
       \end{array}
       \right)_{.}
       \]
By adding the 5th column to the 1st column, we have
\[ T= \left( \begin{array}{ccccc} 
       0 & 0 & 1 & 0 & 0 \\  
       0 & 0 & 0 & 1 & 0 \\
       1 & 0 & 0 & 0 & 1 \\
       d & 0 & 0 & 0 & d \\
       * & *  & *  & * & 0 
       \end{array}
       \right)
       :
       \left(
       \begin{array}{ccccc} 
        1 & 0 & 0 & 0 & 0 \\  
        0 & 1 & 0  & 0 & 0 \\
       0 & 0 &  1 & 0 & 0 \\
       0 & 0 &  0&  1 & 0 \\
       0 & 0 & 0 & 0& 0 
       \end{array}
       \right)_{.}
       \]
   Adding the vector $(-d,1,0,0,0)$ to the 4th row, we have
   \[ T= \left( \begin{array}{ccccc} 
       0 & 0 & 1 & 0 & 0 \\  
       0 & 0 & 0 & 1 & 0 \\
       1 & 0 & 0 & 0 & 1 \\
       0 & 1 & 0 & 0 & d \\
       * & *  & *  & * & 0 
       \end{array}
       \right)
       :
       \left(
       \begin{array}{ccccc} 
        1 & 0 & 0 & 0 & 0 \\  
        0 & 1 & 0  & 0 & 0 \\
       0 & 0 &  1 & 0 & 0 \\
       0 & 0 &  0&  1 & 0 \\
       0 & 0 & 0 & 0& 0 
       \end{array}
       \right)_{.}
       \]
   Thus we have $r(T) \leq 1+4+2=7$.
   If $c=0$, multiplication of a constant to the 5th column, we have
   \[ T= \left( \begin{array}{ccccc} 
       0 & 0 & 1 & 0 & 0 \\  
       0 & 0 & 0 & 1 & 0 \\
       0 & 0 & 0 & 0 & 0 \\
       0 & 0 & 0 & 0 & 1 \\
       * & *  & *  & * & 0 
       \end{array}
       \right)
       :
       \left(
       \begin{array}{ccccc} 
        1 & 0 & 0 & 0 & 0 \\  
        0 & 1 & 0  & 0 & 0 \\
       0 & 0 &  1 & 0 & 0 \\
       0 & 0 &  0&  1 & 0 \\
       0 & 0 & 0 & 0& 0 
       \end{array}
       \right)_{.}
       \]
   Adding the 5th column to the 2nd column we have
   \[ T= \left( \begin{array}{ccccc} 
       0 & 0 & 1 & 0 & 0 \\  
       0 & 0 & 0 & 1 & 0 \\
       0 & 0 & 0 & 0 & 0 \\
       0 & 1 & 0 & 0 & 1 \\
       * & *  & *  & * & 0 
       \end{array}
       \right)
       :
       \left(
       \begin{array}{ccccc} 
        1 & 0 & 0 & 0 & 0 \\  
        0 & 1 & 0  & 0 & 0 \\
       0 & 0 &  1 & 0 & 0 \\
       0 & 0 &  0&  1 & 0 \\
       0 & 0 & 0 & 0& 0 
       \end{array}
       \right)_{.}
       \]
  Adding 1 to the $(3,1)$ cell of $A_{2}$ 
 \[ T= \left( \begin{array}{ccccc} 
       0 & 0 & 1 & 0 & 0 \\  
       0 & 0 & 0 & 1 & 0 \\
       1 & 0 & 0 & 0 & 0 \\
       0 & 1 & 0 & 0 & 1 \\
       * & *  & *  & * & 0 
       \end{array}
       \right)
       :
       \left(
       \begin{array}{ccccc} 
        1 & 0 & 0 & 0 & 0 \\  
        0 & 1 & 0  & 0 & 0 \\
       0 & 0 &  1 & 0 & 0 \\
       0 & 0 &  0&  1 & 0 \\
       0 & 0 & 0 & 0& 0 
       \end{array}
       \right)_{.}
       \]
   Thus we have that $r(T) \leq 1+4+2=7.$ These complete the proof of the theorem.
\end{proof}   
\begin{remark}
The proof technique of  many propositions in this section is so elementary, and there is a possibility that they might have been appeared somewhere. However as far as the authors know, at least, the proof technique for Proposition 2.5 seems to be new. 
\end{remark}
\section{Main Theorem}
In this section we will give an simple inductive proof for the formula that $r(2,n,n) \leq [3n/2]$. First we treat the non-singular case.
\subsection{Estimation by Symmetrization }
In this subsection we show that for the case with  non-singular components  the upper bound $[3n/2]$ for $r(2,n,n)$ is easily proved by the symmetrization method.  
\begin{thm}
If there is contained at least one non-singular matrix in $\langle A,B\rangle$, $r(A;B) \leq [3n/2]$
\end{thm} 

For the proof we prepare two lemmas. 
\begin{lemma} For a $n \times n$ square real matrix  $F$,  there is a factorization of $F=AB^{-1}$ or $F=B^{-1}A$, where $A,B$ are appropriate real symmetric matrices.  
\end{lemma}
\begin{proof} For the proof see Bosch~\cite{Bosch}.
\end{proof}
\begin{lemma}
For a pair of  symmetric matrices $A$ and $B$, if at least one of them is positive definite, they are diagonalizable simultaneously by congruence. That is, there is a matrix P such that  $P^{t}AP=D_{1},P^{t}BP=D_{2}, D_{1}, D_{2}$ are both diagonal matrices.
\end{lemma} 
\begin{proof}
The proof is easy and omitted.
\end{proof}

\begin{proof} 
For a $2 \times n \times n$ tensor $T=(A;B)$, without loss of generality, we assume that  $B$ is non-singular.
By singular value decomposition, multiplying non-singular matrix from both sides, we have that $T=(A;E_{n})$.
Here note that $A$ is transformed by the same operation, without confusion, we use the same symbol $A$.
By using Lemma 1, $A$ is expressed as $A=PQ^{-1}$ where $P$ and $Q$ are appropriate real symmetric  matrices. Therefore, $T=(PQ^{-1};E_{n})$.
From this $T$ is equivalent with $T=(P;Q)$. Since $Q$ is symmetric, it is diagonalized by an orthogonal matrix, 
\[ B=Q=\left(\begin{array}{cccc}
                 \lambda_{1} & 0 & \cdots & 0\\
                     0     &   \lambda_{2} & \cdots & 0 \\ 
                     \cdots & \cdots & \cdots & \cdots \\
                       0 &  0  & 0  & \lambda_{n} 
                       \end{array}
                       \right)
\]
If necessary, multiplying (-1) to both matrices of the tensor, 
at least $n/2$ diagonal elements of $B$ can be assumed as positive. Therefore, 
by adding at most $[n/2]$ positive diagonal elements,  all  diagonal element become positive, and the matrix $B$ becomes  positive definite. Then by Lemma 2 both $A$ and $B$ are diagonalizable simultaneously.  From these the rank of $T$ is less than or equal to $n+[n/2]=[3n/2]$. This completes the proof of Theorem 1.
\end{proof}

\subsection{$r(2,n,n)=[3n/2]$ for even $n$}
In this subsection we will show that for even $n$ the formula is automatically proved by a very simple induction.  Here we prove this briefly. 
First we prove the following lemma.
\begin{lemma}
 Under the assumption that $r(n-2,n-2) \leq [3(n-2)/2]$, it holds that $r(2,n,n+2)=[3n/2]+1$
 \end{lemma}
\begin{proof}
We prove the lemma by induction and so we assume that $r(2,n-2,n-2) \leq [3(n-2)/2].$ Here we 
consider a tensor as $n$ slices of $2 \times n$ matrix. Thus, all symbols denote $n$-dimensional vectors.
We can start from 
\[\left(
\begin{array}{cccccccc} 
\bm{a}_1 & \bm{a}_2 & \cdot & \cdots &  \bm{a}_{n-1} & \bm{a}_n & \bm{0} & \bm{0} \\   
0 & 0 & \bm{b}_{n} & \bm{b}_{n-1} & \cdots & \cdot & \bm{b}_{2} & \bm{b}_{1} 
\end{array}
\right), 
\]
where $\bm{a}_1,,,,\bm{a}_n$ are independent $n$-dimensional vectors and also $\bm{b}_{1},,,,\bm{b}_{n}$ are independent $n$-dimensional vectors. 
Since both of $V=\langle \bm{a}_1,\bm{a}_2 \rangle$ and $W=\langle \bm{b}_{1},\bm{b}_{2} \rangle $ are  $2$-dimensional vector spaces, there is a common $(n-2)$-dimensional vector space $Z$
such that $V \oplus Z=W\oplus Z=R^{n}$. Thus without loss of generality we can write 
\[\left(
\begin{array}{cccccccc} 
\bm{a}_1 & \bm{a}_2 & \bm{z}_{11} & \cdots &  \bm{z}_{1,n-3} & \bm{z}_{1,n-2} & \bm{0} & \bm{0} \\   
\bm{0} & \bm{0} & \bm{z}_{2,n-2} & \bm{z}_{2,n-3} & \cdots & \bm{z}_{2,1} & \bm{b}_{1} & \bm{b}_{2}
\end{array}
\right),
\]
where $\bm{z}_{ij} \in Z$. Hence we have $r(T) \leq r(T_{1})+4,$ where
 \[T_{1}=\left( \begin{array}{cccccccc}
\bm{0} & \bm{0} & \bm{z}_{11} & \cdots &  \bm{z}_{1,n-3} & \bm{z}_{1,n-2} & \bm{0} & \bm{0} \\   
\bm{0} & \bm{0} & \bm{z}_{2,n-2} & \bm{z}_{2,n-3} & \cdots & \bm{z}_{2,1} & \bm{0} & \bm{0}
\end{array}
\right)_{.}
\]
 Since $Z$ is  a $(n-2)$-dimensional  subvector space of $R^{n-1}$ there is a nonsingular matrix $G$ such that
$G\bm{z}_{ij}=(*, *, \cdots, *, 0)^{T}$. Hence
$$ r(T) \leq r(T_{1})+4 \leq r(2,n-2,n-2)+4=[3(n-2)/2]+4=[3n/2]+1,$$
which completes the proof of the statement of the lemma.\\
\end{proof}
Now we begin to prove the following theorem.
\begin{thm}
$r(2,n,n) \leq [3n/2]$ for even $n$. 
\end{thm}

\begin{proof}
As an inductive assumption we assume that $r(2,m,m) \leq [3m/2]$ for all even  $m$ less than $n$. Note that this assumption is assumed through this section.  
If one of $A_{1}$ or $A_{2}$ is non singular, we have already proved the statement of the theorem by the symmetrization method. 
So we assume, both ranks of $A_{1}$ and $A_{2}$ are  singular. Further if one of the ranks is less than $n-2$,  from the previous lemma, we have 
$r(T) \leq r(n-2,n)+2 \leq [3(n-2)/2]+2=[3n/2]-1$. Thus we assume that both of  $A_{1}$ and  $A_{2}$ are of rank $(n-1)$. 

Then we can start from 
\[
\left( \begin{array}{cc} 
       A_{n-1,n-1} & 0_{n-1,1} \\  
        0_{1,n-1} &    0 \\
       \end{array}
       \right)
       :
       \left(
       \begin{array}{ccc}
       B_{n-2,n-2} & 0_{n-2,1} & 0_{n-2,1}  \\
       0_{1,n-2}& 0 & 1 \\  
       0_{1,n-2} & 1 & 0 
       \end{array}
       \right)_{.}
       \]
  where $A_{n^1,n-1}$ is nonsingular . Hence we have
  $r(T) \leq r(T_{1})+2$ where $T_{1}=(A_{n-1,n-1},B_{n-1,n-1})$ with $B_{n-1,n-1}$ below
      \[  \left(
       \begin{array}{cc}
       B_{n-2,n-2} & 0_{n-2,1} \\  
       0_{1,n-2}& 0 \\  
       \end{array}
       \right)_{.}
       \]
 Since $A_{n-1,n-1}$ is nonsingular, $r(T_{1})=[3(n-1)/2]$. And so,
 $$
  r(T) \leq r(T_{1})+2=[3(n-1)/2]+2=[3(2k-1)/2]+2=3k=[3n/2]
$$
This completes the proof of the formula for even $n$. 
\end{proof}
Thus we only need to give a proof for odd $n$. This is very subtle problem to solve.  Therefore we must depart form this simple induction method and goes to the proof based on the following lemma which are also proved by induction. It should be noted that the proof  is applicable both for odd and even $n$.  
 \\ First we need the following lemma for the proof of the main theorem.
\begin{thm}
$$
r(n-1,n) \leq [3n/2]-1
$$
\end{thm}

\begin{proof}
Here we consider a tensor as $(n-1)$ slices of $2 \times n$ matrix. Thus, all symbols denote $(n-1)$-dimensional vectors.
We can start from 
\[\left(
\begin{array}{cccccc} 
\bm{a}_1 & \bm{a}_2 & \cdot & \cdots &  \bm{a}_{n-1}  & \bm{0} \\   
\bm{0} & \bm{b}_{n-1}  & \cdots & \cdot & \bm{b}_{2} & \bm{b}_{1} 
\end{array}
\right),
\]
where $\bm{a}_1,,,,\bm{a}_{n-1}$ are independent $(n-1)$-dimensional vectors and also $\bm{b}_{1},,,,\bm{b}_{n-1}$ are independent $(n-1)$-dimensional vectors. 
Since both of $V=\langle \bm{a}_1\rangle$ and $W=\langle \bm{b}_{1} \rangle$ are  $1$-dimensional vector subspaces of $R^{n-1}$ there is a common $(n-2)$- dimensional vector sub space $Z$ such that $V \oplus Z=W\oplus Z=R^{n-1}$. Thus without loss of generality we can write 
\[\left(
\begin{array}{cccccc} 
\bm{a}_1 &  \bm{z}_{11} & \cdots &  \bm{z}_{1,n-3} & \bm{z}_{1,n-2} & \bm{0}  \\   
\bm{0} &  \bm{z}_{2,n-2} & \bm{z}_{2,n-3} & \cdots & \bm{z}_{2,1} & \bm{b}_{1} 
\end{array}
\right)
\]
where $\bm{z}_{ij} \in Z.$ Hence we have $ r(T) \leq r(T_{1})+2,$ where 
\[ T_{1}=\left( \begin{array}{cccccc}
\bm{0} &  \bm{z}_{11} & \cdots &  \bm{z}_{1,n-3} & \bm{z}_{1,n-2} & \bm{0} \\   
\bm{0} &  \bm{z}_{2,n-2} & \bm{z}_{2,n-3} & \cdots & \bm{z}_{2,1} & \bm{0} 
\end{array}
\right)_{.}
\]
Since $Z$ is  a $(n-2)$-dimensional subvector space of $R^{n-1}$ there is a nonsingular matrix $G$ such that
$G\bm{z}_{ij}=(*,*,\cdots,*,0)^{T}$. Hence
$$
r(n-1,n) \leq r(n-2,n-2)+2=[3(n-2)/2]+2=[3n/2]-1, 
$$
which completes the proof of the lemma. Now we proceed to the proof of the main theorem.
\end{proof}

\begin{thm}
$$r(2,n,n) \leq [3n/2]$$
\end{thm}
\begin{proof}
Let $T=(A_{1}:A_{2})$. We assume $A_{1}$ and $A_{2}$ are of rank $(n-1)$.  Then we can start from 
\[
\left( \begin{array}{cc} 
       A_{n-1,n-1} & 0_{n-1,1} \\  
        0_{1,n-1} &    0 \\
       \end{array}
       \right)
       :
       \left(
       \begin{array}{ccc}
       B_{n-2,n-2} & 0_{n-2,1} & 0_{n-2,1}  \\
       0_{1,n-2}& 0 & 1 \\  
       0_{1,n-2} & 1 & 0 
       \end{array}
       \right)_{.}
       \]
       From this, we have
       $$r(T) \leq r(2,n-1,n)+1$$
 From the previous lemma 
      $$r(T) \leq r(2,n-1,n)+1 \leq [3n/2]-1+1=[3n/2],$$
      which completes the proof of the main theorem.
\end{proof}
\begin{remark} It is known that the reverse inequality holds for some tensors in $T(2,n,n)$, and  in fact it holds that $r(2,n,n)=[3n/2]$.   
\end{remark}

\section{ A generalization to $ 2 \times m\times n $}
In this section we generalize the result in the previous section. The proof  is on the same line.
\begin{thm}
For $m\leq n\leq 2m$ it holds
$r(2,m,n)=m+\lfloor \frac{n}{2}\rfloor$.
\end{thm}
\begin{proof}
It has already known that for some tensor it's rank is greater than or 
equal to $m+\lfloor \frac{n}{2}\rfloor$.
So, we must show $r(2,m,n)\leq m+\lfloor \frac{n}{2}\rfloor$.
If $n\geq 2m$ it is also know that $r(2,m,n)=2m$.
Thus we may assume that $m\leq n<2m$.
We will show by induction on $m$.
Assume that it holds $r(2,k,n)=k+\lfloor \frac{n}{2}\rfloor$
for arbitrary $k<m$ and $k\leq n\leq 2k$.
Consider a $2\times m\times n$ tensor $T$ as $m$ slices of $2\times n$ matrices:
$$
T=\begin{pmatrix}
\bm{x}_1&\bm{x}_2&\cdots&\bm{x}_n\\
\bm{y}_1&\bm{y}_2&\cdots&\bm{y}_n\\
\end{pmatrix}
$$
In the previous section we proved for $m=n$ and 
\par
now we let $m<n<2m$.
We can transform $T$ to
$$
\begin{pmatrix}
\bm{a}_1&\bm{a}_2&\cdots&\bm{a}_s&\bm{0}&\cdots&\bm{0}&\bm{0}&\cdots&\bm{0}\\
\bm{b}_1&\bm{b}_2&\cdots&\bm{b}_s&\bm{b}_{s+1}&\cdots&\bm{b}_{s+t}&\bm{0}&\cdots&\bm{0}\\
\end{pmatrix}
$$
for some $s\leq m$ such that $\bm{a}_1,\ldots,\bm{a}_s$ and $\bm{b}_{s+1},\ldots,\bm{b}_{s+t}$ 
are linearly independent respectively.
Thus the rank of this tensor has an upper bound $r(2,s+t,m)$
which is less or equal to $r(2,n,m)$.
So, we can assume $s+t=n$.
If $\langle \bm{b}_1 \ldots, \bm{b}_s\rangle$ is a subspace of 
$\langle \bm{b}_{s+1},\ldots, \bm{b}_{n} \rangle$, we can transform it to
$$
\begin{pmatrix}
\bm{a}_1&\cdots&\bm{a}_u&\bm{a}_{u+1}&\cdots&\bm{a}_s&\bm{0}&\cdots&\bm{0}\\
\bm{0}&\cdots&\bm{0}&\bm{0}&\cdots&\bm{0}&\bm{b}_{s+1}&\cdots&\bm{b}_{n}\\
\end{pmatrix}
$$
and thus $r(T)\leq n$.
Let suppose that $\langle \bm{b}_1 \ldots, \bm{b}_s\rangle$ is not a subspace of
$\langle \bm{b}_{s+1},\ldots, \bm{b}_{n}\rangle$.
Then we transform it to
$$
\begin{pmatrix}
\bm{a}_1&\cdots&\bm{a}_u&\bm{a}_{u+1}&\cdots&\bm{a}_s&\bm{0}&\cdots&\bm{0}\\
\bm{0}&\cdots&\bm{0}&\bm{b}_{u+1}&\cdots&\bm{b}_s&\bm{b}_{s+1}&\cdots&\bm{b}_{n}\\
\end{pmatrix}
$$
for some $u\leq s-1$ such that $\bm{b}_{u+1},\ldots, \bm{b}_s$ are linearly independent.
Note that $u< s\leq m$ and $n-u \leq m$.
Let $d=\min(u,n-s)$.
Take a vector space $Z$ which has a minimal dimension 
among $Z$ satisfying that
$$
Z+\langle \bm{a}_1,\ldots,\bm{a}_d \rangle=Z+ \langle \bm{b}_{n-d+1},\ldots,\bm{b}_n\rangle
=\langle \bm{a}_1,\ldots,\bm{a}_s,\bm{b}_{u+1},\ldots,\bm{b}_n\rangle.
$$
Then we can transform it to the above form with
$\bm{a}_{d+1}, \ldots, \bm{a}_s, \bm{b}_{u+1}\ldots, \bm{b}_{n-d} \in Z$.
Thus $r(T)$ is less than or equal to $2d+r(2,\dim(Z),n-2d)$.
%%%
Since $\dim(Z)\leq m-d$, by the assumption of the induction,
we have 
\begin{equation*}
\begin{split}
r(T) & \leq 2d+\left(\dim(Z)+\lfloor \frac{n-2d}{2}\rfloor \right) \\
 & \leq m+\lfloor \frac{n}{2}\rfloor.
\end{split}
\end{equation*}
%%%%
We completes the proof when $m<n<2m$.
\end{proof}

\end{document}